\documentclass[a4paper, 11pt]{amsart}

\usepackage{amsmath,amsthm, amscd, amssymb, amsfonts}
\usepackage{latexsym}
\usepackage{graphicx}
\usepackage{color}
\usepackage{indentfirst}
\usepackage[all]{xy}

\usepackage[active]{srcltx}

\theoremstyle{plain}

\newtheorem{theorem}{Theorem}[section]
\newtheorem{lemma}[theorem]{Lemma}
\newtheorem{proposition}[theorem]{Proposition}
\newtheorem{corollary}[theorem]{Corollary}
\newtheorem{conjecture}{Conjecture}

\theoremstyle{definition}
\newtheorem{examples}[theorem]{Examples}

\newtheorem{definition}[theorem]{Definition}

\newcounter{quest}

\newtheorem{question}[quest]{Question}

\theoremstyle{remark}
\newtheorem{remark}[theorem]{Remark}

\newcommand{\D}{{\mathcal D}}
\newcommand{\mh}{{\mathcal M}^H}
\newcommand{\ma}{{\mathcal M}^A}
\newcommand{\mbr}{{\mathcal M}^B}
\newcommand{\mcr}{{\mathcal M}^C}
\newcommand{\mb}{{}^B\hspace{-1pt}{\mathcal M}}
\newcommand{\mc}{{}^C\hspace{-1pt}{\mathcal M}}
\newcommand{\mk}{{\mathcal M}^K}
\newcommand{\mkfd}{{}_K\hspace{-1pt}{\mathcal M}\hspace{-2pt}_{f}}
\newcommand{\mhfd}{{}_H\hspace{-1pt}{\mathcal M}\hspace{-2pt}_{f}}
\newcommand{\hm}{{}^H\hspace{-1pt}{\mathcal M}}
\newcommand{\km}{{}^K\hspace{-1pt}{\mathcal M}}
\newcommand{\hcerom}{{\mathcal M}^{H^0}}
\newcommand{\kcerom}{{\mathcal M}^{K^0}}

\newcommand{\Ss}{{\mathcal S}}
\newcommand{\C}{{\mathcal C}}
\newcommand{\Oc}{{\mathcal O}}

\def\Oe{\Oc_{\epsilon}(G)}

\newcommand{\Na}{\mathbb{N}}

\newcommand{\Rat}{\operatorname{Rat}}
\newcommand{\Rad}{\operatorname{Rad}}

\newcommand{\ku}{\Bbbk}

\newcommand\toba{{\mathfrak B }}

\newcommand{\YDL}{{}_{H_{[0]}}^{H_{[0]}}{\mathcal YD}}
\newcommand{\ydh}{{}_{H_{0}}^{H_{0}}{\mathcal YD}}
\newcommand{\ydl}{{}_{L}^{L}{\mathcal YD}}

\newcommand{\gr}{\operatorname{gr}}
\newcommand{\ord}{\operatorname{ord}}
\newcommand\co{\operatorname{co}}
\newcommand\id{\operatorname{id}}

\newcommand{\Hom}{\operatorname{Hom}}
\newcommand{\Gr}{\mathfrak{gr}}
\newcommand{\Ind}{\mathsf{Ind}\,}
\newcommand{\Res}{\mathsf{Res} \,}

\def\pf{\begin{proof}}

\def\epf{\end{proof}}

\definecolor{rojo}{rgb}{1,0,0}
\definecolor{verde}{rgb}{0,1,0}

\begin{document}
\renewcommand{\baselinestretch}{1.2}

\thispagestyle{empty}

\title[The structure of (co-frobenius) Hopf algebras]
{On the structure of (co-frobenius) Hopf algebras}

\author[Andruskiewitsch, Cuadra]{Nicol\'as Andruskiewitsch, Juan Cuadra}

\address{\noindent N. A.: Facultad de Matem\'atica, Astronom\'{\i}a y F\'{\i}sica,
Universidad Nacional de C\'ordoba. CIEM -- CONICET. 
Medina Allende s/n (5000) Ciudad Universitaria, C\'ordoba,
Argentina}
\email{andrus@famaf.unc.edu.ar}

\address{J. C.: Universidad de Almer\'\i a, Dpto. \'Algebra y An\'alisis Matem\'atico.
E04120 Almer\'\i a, Spain}

\email{jcdiaz@ual.es}


\date{\today}

\begin{abstract}
We introduce a new filtration on Hopf algebras, the standard filtration, generalizing the coradical filtration. Its zeroth  term, called the Hopf coradical, is the subalgebra generated by the coradical. We give a structure theorem: any Hopf algebra with injective antipode is a deformation of the bosonization of the Hopf coradical by its diagram, a connected graded Hopf algebra in the category of Yetter-Drinfeld modules over the latter. We discuss the steps needed to classify Hopf algebras in suitable classes accordingly. For the class of co-Frobenius Hopf algebras, we prove that a Hopf algebra is co-Frobenius if and only if its Hopf coradical is so and the diagram is finite dimensional. We also prove that the standard filtration of such Hopf algebras is finite. Finally, we show that extensions of co-Frobenius (resp. cosemisimple) Hopf algebras are
co-Frobenius (resp. cosemisimple).
\end{abstract}

\maketitle

\section*{Introduction}

There are few general techniques to deal with the classification of Hopf algebras; one of them is the so-called Lifting Method \cite{AS-cambr}
under the assumption that the coradical is a subalgebra. In the present paper we propose to extend this technique by considering the subalgebra generated by the coradical, called the Hopf coradical, and the related standard filtration, that is a generalization of the coradical filtration. Its zeroth term is the Hopf coradical while the remaining ones are iterative wedge operations of it. The standard filtration of a Hopf algebra $H$ is always a Hopf algebra filtration, provided that the antipode is injective, and we can consider its associated graded Hopf algebra $\gr H$. The latter is a bosonization of the Hopf coradical $H_{[0]}$ by a connected graded Hopf algebra $R$ in the category of Yetter-Drinfeld $H_{[0]}$-modules (the diagram of $H$). Then $H$ is a deformation or quantization of $\gr H$ for a suitable cohomology theory. We summarize our considerations in Theorem \ref{th:structure}, that can be thought of either as a structure theorem for Hopf algebras with injective antipode, or as a proposal for the classification of Hopf algebras in suitable classes (e.g., those of finite dimensional, or co-Frobenius, or finite Gelfand-Kirillov dimension Hopf algebras). We discuss the different problems to be solved for the success of this proposal in Section \ref{sec:standard}. \par \smallskip

In Section \ref{sec:cofrob} we focus on the class of co-Frobenius Hopf algebras with the techniques just introduced. These are Hopf algebras having  nonzero integral and there exist relevant examples of them: either finite dimensional or cosemisimple Hopf algebras, coordinate algebras of certain algebraic groups \cite{Su2} or more generally group schemes \cite{Don2}, families of quantum groups at a root of one \cite{APW,AD}, and quantum groups attached to some Hecke symetries \cite{Hai}. This notion can be rephrased in representation theoretic terms: a  Hopf algebra is co-Frobenius iff the injective hulls of the simple comodules are finite dimensional iff the projective cover of any comodule do exist. \par \smallskip

Our second main result, Theorem \ref{teo:cofrob}, particularizes the previous structure theorem as follows: a Hopf algebra is co-Frobenius if and only if its Hopf coradical is so and the diagram is finite dimensional. This reduces the classification in this class to those which are generated by cosemisimple coalgebras and the finite dimensional graded connected Hopf algebras in the corresponding categories of Yetter-Drinfeld modules. It is natural to investigate cosemisimple Hopf algebras as part of this program. Towards this end, we study extensions of co-Frobenius Hopf algebras.
Given an extension of Hopf algebras $1 \rightarrow A \rightarrow B \rightarrow C \rightarrow 1$ with $B$ faithfully coflat as a $C$-comodule,
Theorem \ref{exseq} asserts that $B$ is co-Frobenius if and only if $A$ and $C$ are co-Frobenius. This result is used to detect all co-Frobenius quotients of quantized coordinate algebras of simple algebraic groups at root of one, Examples \ref{examples}. In our last main result, Theorem \ref{exseq2}, we prove that $B$ has a nonzero left integral that restricted to $A$ is nonzero if and only if $A$ is co-Frobenius and $C$ is cosemisimple. We derive from this that $B$ is cosemisimple if and only if $A$ and $C$ so are. New characterizations of co-Frobenius Hopf algebras are established to achieve these results, Theorems \ref{newchar1} and \ref{newchar2}. \par \smallskip

A conjecture posed in \cite{AD} states that a co-Frobenius Hopf algebra has finite coradical filtration. Related to this problem, in Theorem \ref{teo:cofrob} we also show that a Hopf algebra is co-Frobenius if and only if the Hopf coradical is so and its standard filtration is finite. After acceptance of the present paper, it was proved in \cite[Theorem 1.2]{ACE} that the conjecture is true.  \par \smallskip

Contents of Section 1 were presented by N. Andruskiewitsch at the meetings ``Conference in Hopf algebras and Noncommutative Algebra'', Sde-Boker (Israel), May 24-27, 2010 and ''XXI Escola de Algebra'', Brasilia (Brazil), July 25-31, 2010. Main results of Section 2 were expounded by J. Cuadra at the conference ``Quantum groups: Galois and integration techniques'', Clermont-Ferrand (France), August 30-September 3, 2010. \par \smallskip

\subsection*{Conventions and notations}  Our main references for the theory of Hopf algebras are \cite{Mon, Sw}. We shall work over a ground field $\ku$. Let $C$ be a coalgebra with comultiplication $\Delta$ and counit $\varepsilon$. For subspaces $D, E \subset C$ recall from \cite[Proposition 9.0.0]{Sw} that the {\it wedge} of $D$ and $E$ is defined to be $D \wedge E = \{c\in C: \Delta(c) \in D\otimes C + C \otimes E\}$. Using the dual algebra $C^*$, it is $D \wedge E = (D^{\perp}E^{\perp })^{\perp}$, where $\perp$ stands for the annihilator subspace (in $C^*$ and $C$). We inductively define $\wedge^{0} D = D$ and $\wedge^{n+1} D = (\wedge^{n} D) \wedge D$ for $n>0.$ \par \smallskip

We denote by $\Ss$ the antipode of any Hopf algebra $H$ and by $H^+$ the kernel of the counit. We shall use that $\Ss(D \wedge E) \subseteq \Ss(E) \wedge \Ss(D)$ and that this is an equality when $\Ss$ is bijective.

\section{The standard filtration}\label{sec:standard}
Let $H$ be a Hopf algebra.
We shall consider several invariants of $H$. The first one, already
present in \cite{Sw}, is the {\em coradical filtration} $\{H_{n}\}_{n\geq 0}$, whose terms are defined as follows:
\begin{itemize}
  \item $H_0$ is the coradical, i.e., the sum of all simple subcoalgebras of $H$.
  \item $H_{n}=\wedge^{n+1} H_0$.
\end{itemize}
These are coalgebra versions of the Jacobson radical and its powers; indeed, $H_0=J^{\perp}$ and $H_{n}=(J^{n+1})^{\perp}$, where $J$ denotes the Jacobson radical of $H^*$. The coradical filtration is a coalgebra filtration. Furthermore, if $H_0$ is a Hopf subalgebra, then it is also an algebra filtration \cite[Lemma 5.2.8]{Mon}, and its associated graded coalgebra $\gr H=\oplus_{n\geq 0} H_{n}/H_{n-1}$ is a graded Hopf algebra ($H_{-1}=0$).
Let $\pi:\gr H \rightarrow H_{0}$ be the homogeneous projection; since it splits the inclusion of $H_{0}$ in $\gr H$, the diagram $$R = (\gr H)^{\co \pi} = \{x\in \gr H: (\id \otimes \pi)\Delta (x) = x\otimes 1\}$$
turns out to be a Hopf algebra in the category $\ydh$ of Yetter-Drinfel'd $H_{0}$-modules and
$\gr H\cong R \# H_{0}$. Here $\#$ stands for the Radford biproduct or bosonization, see for example \cite{AS-cambr}. The study of the diagram is central for the understanding of Hopf algebras whose coradical is a Hopf subalgebra. But this is not always the case, and the main goal of this paper is to propose a new approach in the general situation. We start by defining a new filtration, the \emph{standard filtration} $\{H_{[n]}\}_{n \geq 0}$, as follows:
\begin{itemize}
  \item The \emph{Hopf coradical} $H_{[0]}$ is the subalgebra generated by $H_0$.
  \item $H_{[n]}=\wedge^{n+1} H_{[0]}.$
\end{itemize}
By convenience, we set $H_{[-1]}=0$. Of course, $H_{[0]} = H_0$ just means that the latter is a subalgebra;
then it is a Hopf subalgebra and the coradical filtration coincides with the standard one. The basic properties of the standard filtration are collected in the next result.

We assume throughout this section that $\Ss(H_0) \subseteq H_0$; this holds, for instance, if $\Ss$ is injective. Actually, we are mostly interested in Hopf algebras with bijective antipode.

\begin{lemma}\label{lema:filt} With notation as above:
\begin{enumerate}
\item[(i)] $H_{[0]}$ is a Hopf subalgebra of $H$ and its coradical is $H_0$.
\item[(ii)] $H_n \subseteq H_{[n]}$ and $\{H_{[n]}\}_{n \geq 0}$ is a Hopf algebra filtration
of $H$.
\item[(iii)] If $\Ss$ is bijective, then $\Ss(H_{[n]})=H_{[n]}.$
\end{enumerate}
\end{lemma}

\pf (i) We know that $H_{[0]}=\bigcup_{r \geq 0} H_0^{(r)},$ where $H_0^{(r)}=H_0 \stackrel{r_{}}{...} H_0$
for $r > 0$ and $H_0^{(0)}=\ku$. Then $H_{[0]}$ is a subcoalgebra of $H$ because each $H_0^{(r)}$ is so.
Since $\Ss(H_0) \subseteq H_0$ by assumption,  $\Ss(H_0^{(r)}) \subseteq H_0^{(r)}$ and thus $\Ss(H_{[0]}) \subseteq H_{[0]}$. For the second statement, the coradical of $H_{[0]}$ is $H_{[0]} \cap H_0=H_0.$ \par \medskip

(ii) This can be similarly proved as \cite[Lemma 5.2.8]{Mon}; we include the proof for the sake of completeness.
Each $H_{[n]}$ is a subcoalgebra of $H$, because it is defined as an iterative wedge of subcoalgebras, and
$H_{[n]} \subseteq H_{[n+1]}$, \cite[Proposition 9.0.0(i)]{Sw}. Moreover,
from $H_n =\wedge^{n+1} H_0 \subseteq \wedge^{n+1} H_{[0]}=H_{[n]}$ and $H=\bigcup_{n \geq 0} H_n$ we obtain $H=\bigcup_{n \geq 0} H_{[n]}$.
Since $H_{[n]}=\wedge^{n+1} H_{[0]},$ by \cite[Theorem 9.1.6]{Sw}, $\Delta(H_{[n]}) \subseteq \sum_{i=0}^{n} H_{[i]} \otimes H_{[n-i]},$ showing that $\{H_{[n]}\}_{n \geq 0}$ is a coalgebra filtration. We now prove that it is an algebra filtration, that is, $H_{[n]}H_{[m]} \subseteq H_{[n+m]}$ for all $n,m \geq 0$. For $n=0$, it follows by induction on $m$ and the following computation:
$$\begin{array}{ll}
\Delta(H_{[0]}H_{[m]}) & \subseteq (H_{[0]} \otimes H_{[0]})(H_{[0]} \otimes H_{[m]}+H_{[m]} \otimes H_{[m-1]}) \vspace{3pt} \\
 & \subseteq H_{[0]}H_{[0]} \otimes H_{[0]}H_{[m]}+H_{[0]}H_{[m]} \otimes H_{[0]}H_{[m-1]} \vspace{3pt} \\
 & \subseteq H_{[0]} \otimes H + H \otimes H_{[m-1]}
\end{array}$$
Analogously, $H_{[n]}H_{[0]} \subseteq H_{[n]}$ for all $n \geq 0$. To prove the general statement, we apply induction on $n$ and $m$. A computation similar to the preceding one shows by a recursive argument that $H_{[n]}H_{[m]} \subseteq H_{[n+m-1]} \wedge H_{[0]}=H_{[n+m]}$. Finally, since $\Ss$ is an anti-coalgebra map, by induction we have $$\Ss(H_{[n]})=\Ss(H_{[0]} \wedge H_{[n-1]}) \subseteq \Ss(H_{[n-1]}) \wedge \Ss(H_{[0]}) \subseteq H_{[n-1]} \wedge H_{[0]}=H_{[n]}.$$ \par \medskip

(iii) Use that for $\Ss$ bijective, $\Ss(D \wedge E)=\Ss(E) \wedge \Ss(D)$ for any pair of subspaces $D,E$ of $H$. \epf \par \bigskip

Thanks to the previous lemma, we may consider the {\it graded Hopf algebra $\gr H=\oplus_{n \geq 0} H_{[n]}/H_{[n-1]}$ associated with the standard filtration}. As before, if $\pi:\gr H \rightarrow H_{[0]}$ is the homogeneous projection, that splits the inclusion of $H_{[0]}$ in $\gr H$, then the \emph{diagram} $R = (\gr H)^{\co \pi}$ is a Hopf algebra in the category $\YDL$ of Yetter-Drinfel'd $H_{[0]}$-modules and
\begin{equation}\label{eqn:R-bosonization}
\gr H \cong R \# H_{[0]}. \medskip
\end{equation}

For $n \geq 0$ set $\gr^n H= H_{[n]}/H_{[n-1]}$, the homogenous component of degree $n$ in $\gr H$. We are going to see that the filtration of $\gr H$ associated with the grading and the standard filtration of $\gr H$ coincide.

\begin{proposition}\label{stfilgrh}
$(\gr H)_{[n]}=\oplus_{i\leq n}\gr^i H$ for all $n \geq 0$.
\end{proposition}

\pf The proof is similar to that of \cite[Lemma 2.3]{AS-p3}, where this result is established when $H_0$ is a subalgebra. We proceed by induction on $n$. The filtration attached to the grading is a coalgebra filtration. By \cite[Proposition 11.1.1]{Sw}, $(\gr H)_0 \subseteq \gr^0 H=H_{[0]}$. From here, $(\gr H)_{[0]}\subseteq H_{[0]}$. On the other hand, $H_{0}$ is a cosemisimple subcoalgebra of $\gr H$ (as a subcoalgebra of $H_{[0]}$). Hence $H_0 \subseteq (\gr H)_0$ and consequently $\gr^0 H=H_{[0]} \subseteq (\gr H)_{[0]}$. \par \smallskip

The following computation shows that $\oplus_{i\leq n}\gr^i H \subseteq (\gr H)_{[n]}:$
$$\begin{array}{ll}
\Delta(\gr^n H) & \subseteq \oplus_{l=0}^n \gr^l H \otimes \gr^{n-l} H \vspace{3pt} \\
 & \subseteq \gr^0 H \otimes \gr H+\gr H \otimes (\oplus_{i\leq n-1}\gr^i H) \vspace{3pt} \\
 & = (\gr H)_{[0]} \otimes \gr H + \gr H \otimes (\gr H)_{[n-1]}.
\end{array}$$

To prove the other inclusion, we observe that $(\gr H)_{[n]}$ is a graded subspace, that is, $(\gr H)_{[n]} = \oplus_{m\geq 0}(\gr^m H \cap (\gr H)_{[n]})$;
the wedge of two graded subspaces is graded. Thus, it suffices to check that $\gr^m H \cap (\gr H)_{[n]} = 0$ for $m>n$. Pick $0\neq\overline{h} \in H_{[m]}/H_{[m-1]}$ and write $\Delta(h)=x+y+z$ with $x \in \sum_{i=0}^{n-1} H_{[i]} \otimes H_{[m-i]}, y \in H_{[m]} \otimes H_{[0]}$ and $z \in \sum_{i=n}^{m-1} H_{[i]} \otimes H_{[m-i]}$. Applying the corresponding projections defining the comultiplication of $\gr H$ \cite[page 229]{Sw} we can write $\Delta(\overline{h})=\bar{x}+\bar{y}+\bar{z}$ with $\bar{x} \in \oplus_{i=0}^{n-1} \gr^i H \otimes \gr^{m-i} H$, $\bar{y} \in \gr^m H \otimes \gr^0 H$ and $\bar{z} \in \oplus_{i=n}^{m-1} \gr^i H \otimes \gr^{m-i} H$. We claim that $\bar{z} \neq 0$. Otherwise, $z \in
\sum_{i=n}^{m-1} H_{[i]} \otimes H_{[m-1-i]}$ and hence $$\Delta(h) \in H_{[n-1]} \otimes H_{[m]}+H_{[m]} \otimes H_{[m-n-1]}.$$ This yields
$h \in H_{[n-1]} \wedge H_{[m-n-1]}=H_{[m-1]}$ and hence $\bar{h}=0$, a contradiction. Since $\bar{z} \neq 0$, we get that $\Delta(\overline{h}) \notin (\gr H)_{[n-1]}  \otimes \gr H+ \gr H \otimes (\gr H)_{[0]}.$ Therefore, $\bar{h} \notin (\gr H)_{[n]}$. \epf

From the preceding, we deduce that $(\gr H)_n \subseteq \oplus_{i\leq n}\gr^n H$. But the latter is not an equality in general; in other words, $\gr H$ is not coradically graded. \par \medskip

The diagram inherits the grading from $\gr H$, that is, $R= \oplus_{n \geq 0} R^n$ where $R^n= R \cap \gr^n H.$ With respect to this grading, $R$ is a graded Hopf algebra in $\YDL$, $\gr^i H= R^i \# H_{[0]}$ for every $i \geq 0$ and $R_0 = R^0=\ku 1$, see \cite[Lemma 2.1]{AS-p3}. Furthermore,  $R^1 \subseteq P(R)$ as in the proof of \cite[Lemma 2.4]{AS-p3}. \par \bigskip

To sum up this discussion, we have the following structure theorem.

\begin{theorem}\label{th:structure}
Any Hopf algebra with injective antipode is a deformation of the bosonization of a Hopf algebra generated by a cosemisimple coalgebra by a
connected graded Hopf algebra in the category of Yetter-Drinfeld modules over the latter. \qed
\end{theorem}

To provide significance to this result, we should address some fundamental questions.
Suppose that we aim to classify all Hopf algebras $H$ in a class $\C$, that is \emph{suitable} in the following sense:

\begin{enumerate}
\item $H$ belongs to $\C$ iff $\gr H$ belongs to $\C$; \vspace{3pt}
\item If $H$ belongs to $\C$, then $H_{[0]}$ belongs to $\C$.
\end{enumerate}
Typically, the classes of finite dimensional, or with finite Gelfand-Kirillov dimension, or co-Frobenius Hopf algebras are suitable. See Section 2 for the latter class.

\begin{question}
Let $C$ be a cosemisimple coalgebra compatible with the class in the appropriate sense and $S: C\to C$ an injective anti-coalgebra morphism (in the typical examples, one should assume $S$ bijective). Classify all Hopf algebras $L$ generated by $C$, belonging to the class $\C$, and such that $\Ss {\vert_C} = S$.
\end{question}

\begin{question}
Given $L$ as in the previous item, classify all connected graded Hopf algebras $R$ in $\ydl$ such that $R\# L$ belongs to $\C$.
\end{question}

\begin{question}\label{que:deformations}
Given $L$ and $R$ as in previous items, classify all deformations or liftings, that is, classify all Hopf algebras $H$
such that $\gr H \cong R\# L$.
\end{question}

\begin{remark}
There is an alternative dual approach to the one shown before. Namely, let $H$ be a Hopf algebra with surjective antipode and let $J$ denote its Jacobson radical. Let us consider
$$J_{\omega} = \bigcap_{m \geq 0} \wedge^m J.$$
This the largest Hopf ideal contained in $J$. In the finite dimensional case, $J_{\omega}$ was studied in \cite{CH}.
Consider the filtration by Hopf ideals $(J_{\omega}^n)_{n\geq 0}$ and the associated graded Hopf algebra
$\Gr\hspace{2pt} H = \oplus_{n\geq 0} J_{\omega}^n  / J_{\omega}^{n+1}$, where $J_{\omega}^0=H$.
If $H$ is finite dimensional, then $(\Gr\hspace{2pt} H)^* \cong \gr (H^*)$.
However, this setting might be more convenient for the classification of quasi-Hopf algebras, \cite{EG,An}.
\end{remark}

\bigbreak
\subsection{Hopf algebras generated by cosemisimple coalgebras}\label{subsect:gen-by-coalg}
In this subsection, we assume that $\ku$ is an algebraically closed field of characteristic 0. We discuss what is known about the Question I.
Notice that this question contains the classification of all semisimple Hopf algebras, that is largely open, except for some dimensions.
\par \medskip

It is convenient to use the following terminology.

\begin{definition}\label{def:multiplicative} A basis $(e_{ij})_{1\le i,j\le m}$ of a coalgebra $C$ will be called a
\emph{multiplicative} matrix if $\Delta(e_{ij}) = \sum_{p=1}^m e_{ip} \otimes e_{pj}$ and $\varepsilon(e_{ij}) = \delta_{ij}$, the Kronecker symbol.
\end{definition}

We recall now a remarkable result of \c{S}tefan, used in  classification results of low dimensional Hopf algebras \cite{St, Na-12, GV}.

\begin{theorem}\label{teo:stefan} \cite[Theorem 1.5]{St}
Let $H$ be a Hopf algebra and $C$ an $\Ss$-invariant 4-dimensional simple subcoalgebra. If $1 < \ord (\Ss\hspace{-1pt}\vert_{_C}^{\hspace{2pt} 2}) = n < \infty$, then there are a root of unity $\omega$ and a multiplicative matrix $(e_{ij})_{1\le i,j\le 2}$ such that $\ord(w^2) = n$ and $e_{ij}$ satisfy all relations defining $\Oc_{\sqrt{-\omega}}(SL_2(\ku))$. In particular, there is a Hopf algebra morphism $\Oc_{\sqrt{-\omega}}(SL_2(\ku)) \to H$,
which is surjective if $C$ generates $H$ as an algebra.
\qed
\end{theorem}

This raises the question of classifying all quantum subgroups of the quantum group $SL_2$, that is, the quotient Hopf algebras of $\Oc_{q}(SL_2(\ku))$. This problem was considered in \cite{Pod95}. The determination of all quantum subgroups of a quantum group at a root
of one or, in equivalent terms, to determine all Hopf algebra quotients of a quantized coordinate
algebra at a root of one (over $\mathbb C$), was accomplished in \cite{Mul00}, for finite dimensional quotients of
the quantum group $SL_N$, and in \cite{AG}, for quantum versions of simple groups. At the present moment, there is no intrinsic condition describing these quotients, as in the beautiful result of \c{S}tefan for $SL_2$.

\begin{definition}\cite[Lemma 2]{R2} \label{def:H(T)}
Let $C$ be a coalgebra and $S \in GL(C)$ an anti-coalgebra map. The algebra $$\mathcal H(S) := T(C) / \langle c_{(1)}S(c_{(2)}) -\varepsilon (c), S(c_{(1)})c_{(2)} -\varepsilon (c): c\in C\rangle$$ is a Hopf algebra, with comultiplication induced by that of $C$ and antipode induced by $S$, which satisfies the following universal property: if $K$ is a Hopf algebra with antipode ${\mathcal S}_K$ and $f:C \rightarrow K$ is a coalgebra map such that $\mathcal{S}_Kf=fS,$ then there is a unique Hopf algebra map $\tilde{f}:\mathcal{H}(S) \rightarrow K$ such that $\tilde{f} \vert_C=f$.
\end{definition}

Given $s\in \Na$, let $1 < d_1 <\dots < d_s$ and $n_1, \dots, n_s$ be natural numbers. For $1\le r \le s$, let $F_r \in GL_{d_r}(\ku)$.
We consider the coalgebras
$$
D_r =(C_{d_r})^{n_r},\qquad C = \oplus_{r=1}^s D_r,
$$
where $C_{d_r}$ is a comatrix coalgebra of dimension $d_r^2$. Fix $(e^{r,k}_{ij})_{1\le i,j\le d_r}$ a multiplicative matrix of the $k$-th copy of
$C_{d_r}$ in $D_r$, for any $k$, and define $S_r \in GL(D_r)$ by
$$
S_r(e^{r,k}_{ij}) = \begin{cases} e^{r,k + 1}_{ji}, & 1\le k < n_r;
\\ a_{ij}, &  k = n_r.\end{cases}
$$
where $A = (a_{ij})$ is given by $A = F_r (e^{r,1}_{ji}) F_r^{-1}$. Let $S=\oplus_{r=1}^s S_r \in GL(C)$.
We denote $\mathcal H(F_r, n_r)_{1\le r \le s} = \mathcal H(S)$. This definition is a generalization of the one in \cite{B};
a similar construction in the setting of Hopf $C^*$-algebras was introduced in \cite{W}. See also \cite{VDW, BB, BN} for variations and applications.

\begin{question}
Compute the Hopf algebra quotients of $\mathcal H(F_r, n_r)_{1\le r \le s}$ in suitable classes (e.g., finite dimensional, or with finite Gelfand-Kirillov dimension, or co-Frobenius).
\end{question}

\subsection{Connected braided Hopf algebras}\label{subsect:diagram}
We point out here the connection of Question II with Nichols algebras. Let $L$ be a Hopf algebra generated by a cosemisimple coalgebra in our class $\C$ and let $\C_L$ be the class of connected braided Hopf algebras $R$ in $\ydl$ such that $R\# L \in \C$.

The most relevant examples of connected braided Hopf algebras in $\ydl$ are the Nichols algebras: given $V\in \ydl$,
there exists a unique (up to isomorphisms) connected braided Hopf algebra $\toba(V) = \oplus_{n\geq 0}\toba^n(V)$ with the properties
$$
V \cong \toba^1(V)=P(\toba(V)) \hspace{5pt} \text{and $V$ generates $\toba(V)$ as an algebra.}
$$

If $R$ is a connected braided Hopf algebra in $\ydl$, then there is a canonical subquotient Nichols algebra $\toba(V)$, namely
$V = R^1$. Therefore, if the class $\C_L$ is closed under subquotients, then it would be important to solve the following problem.

\begin{question}\label{que:nichols}
Given $L$ as above, classify all Nichols algebras in $\C_L$.
\end{question}

It would be interesting to understand how to construct any braided connected Hopf algebra as a suitable extension of Nichols algebras.
For generalities on extensions in categories of Yetter-Drinfeld modules, see \cite{Bes, BD}.

\subsection{Liftings or deformations}\label{subsect:lifting}
As for Question \ref{que:deformations}, the classification of all Hopf algebras $H$
such that $\gr H \cong R\# L$, with $R$ and $L$ as above, is a particular instance of the general problem of detecting all filtered objects with a fixed graded object $G$. These objects are usually called {\it deformations} or {\it quantizations} of $G$, and they are controlled with a suitable cohomology theory. In the Hopf algebra case, they are called liftings \cite{AS-cambr} and the pertinent cohomology theory is that of \cite{GeS, GeS2}, see \cite{chinos, MaW}.

\section{Co-Frobenius Hopf algebras}\label{sec:cofrob}

Let $H$ be a Hopf algebra. We will denote the category of left (resp. right) $H$-comodules by $\hm$ (resp. $\mh$).
Given $M \in \hm$, throughout this section, $E_H(M)$ stands for the injective hull of $M$.
If $S \in \hm$ is simple, we can always take $E_H(S)$ as a left coideal of $H$, see \cite[1.5g]{G} or \cite[Corollary 2.4.15]{DNR}.
Recall that a {\it left integral} for $H$ is a linear map $\int:H \rightarrow \ku$ such that
$\alpha \cdot \int=\alpha(1_H)\int$ for all $\alpha \in H^*$. Equivalently, $\int(h_{(2)})h_{(1)}= \int(h)1_H$
for all $h \in H$. This is just saying that $\int:H \rightarrow \ku$ is a left $H$-comodule map.
Let $\Rat(H^*)$ denote the maximal rational submodule of $H^*$,  as left $H^*$-module.

\begin{theorem}\label{oldchar}
The following statements are equivalent:
\begin{enumerate}
\item[(i)] $H$ has a nonzero left integral.
\item[(ii)] $\Rat(H^*) \neq 0$.
\item[(iii)] $E_H(S)$ is finite dimensional for every $S \in \hm$ simple.
\item[(iv)] $E_H(\ku)$ is finite dimensional.
\item[(v)] $\hm$ has a nonzero finite dimensional injective object.
\end{enumerate}
\end{theorem}

\pf (i) $\Leftrightarrow$ (ii) $\Leftrightarrow$ (iii) is \cite[Theorem 3]{L}, \emph{cf.} also \cite[2.10]{Sw2}. (iii) $\Rightarrow$ (iv), (iv) $\Rightarrow$ (v) are evident.\footnote{Direct proof of (iv) $\Rightarrow$ (iii): if $S \in \hm$ is simple, then $S \otimes E_H(\ku)$ is injective and contains $S$, so $E_H(S)$ is a subcomodule of $S \otimes E_H(\ku)$ and consequently finite dimensional.}
(v) $\Rightarrow$ (i) is \cite[Proposition 2.3]{DN}. \epf

A Hopf algebra satisfying any of these statements is called {\it co-Frobenius}. Other characterizations may be found in \cite[Theorem 5.3.2]{DNR}; some new ones are given in Theorems \ref{newchar1} and \ref{newchar2}. All these characterizations are equivalent to their right versions, that will be used but not explicitly stated.

\subsection{The standard filtration of co-Frobenius Hopf algebras}

In \cite[Corollary 2]{R} Radford showed that if $H$ is a co-Frobenius Hopf algebra whose coradical $H_0$ is a subalgebra, then $H$ has finite coradical filtration. This is a consequence of his beautiful result:

\begin{theorem}\label{teo:radford-beuatiful}\cite[Proposition 4]{R}
Let $H$ be a co-Frobenius Hopf algebra. Then $H=H_0E_H(\ku)$. \qed
\end{theorem}

As $H_0$ is a subalgebra the coradical filtration $\{H_n\}_{n \geq 0}$ is an algebra filtration. Since $E_H(\ku)$ is finite dimensional, it embeds in $H_m$ for some $m$ and hence $H=H_0E_H(\ku) \subseteq H_0H_m \subseteq H_m$. Note also that $H$ is finitely generated as a left $H_0$-module. In \cite{AD} the relation between co-Frobenius Hopf algebras and the finiteness of the coradical filtration was again analyzed. In that paper, an alternative proof of this fact was provided, it was proved that a Hopf algebra with finite coradical filtration is co-Frobenius \cite[Theorem 2.1]{AD} and the following conjecture was posed:

\begin{conjecture}\cite[page 153]{AD}\label{conj}
The coradical filtration of a co-Frobenius Hopf algebra is finite.\footnote{As of February 2012, it is known that the conjecture is true, see \cite[Theorem 1.2]{ACE}.}
\end{conjecture}

In this subsection we generalize the above-mentioned results by proving that a Hopf algebra $H$ is co-Frobenius if and only if the Hopf coradical $H_{[0]}$ is co-Frobenius and the standard filtration is finite. We will also show that this finiteness condition is reflected in the fact that the diagram $R$ in (\ref{eqn:R-bosonization}) is finite dimensional. In the proof of these results we will need the following new characterization of co-Frobenius Hopf algebras.

\begin{theorem}\label{newchar1}
The following assertions are equivalent:
\begin{enumerate}
\item[(i)] $H$ is co-Frobenius.
\item[(ii)] Every nonzero  $H$-comodule has a nonzero finite dimensional quotient.
\item[(iii)] Every nonzero injective  $H$-comodule has a nonzero finite dimensional quotient.
\item[(iv)] There is an injective  $H$-comodule which has a nonzero finite dimensional quotient.
\end{enumerate}
\end{theorem}

\pf (i) $\Rightarrow$ (ii) Let $0\neq M \in \hm$. Then $E_H(M) \cong \oplus_{i \in I} E_H(S_i)$ where $\{S_i\}_{i \in I}$ is a set of simple subcomodules of $M$ \cite[1.5h]{G}. By hypothesis, $\dim E_H(S_i) < \infty$ for every $i \in I$. Fix $j \in I$ and set $N=\oplus_{i \in I-\{j\}} E_H(S_i)$. Then $0 \neq E_H(M)/N \cong E_H(S_j)$ is finite dimensional, hence $M/M \cap N$ too. We show that $M/M \cap N \neq 0$. If $M \cap N=M$, then $M \cap E_H(S_j) =M \cap N \cap E_H(S_j)=0$, contradicting the fact that $M$ is essential in $E_H(M)$. \par \medskip

(ii) $\Rightarrow$ (iii) and (iii) $\Rightarrow$ (iv) are obvious. \par \medskip

(iv) $\Rightarrow$ (i) Let $M\in \hm$ be such injective comodule and $g:M \rightarrow P$ a surjective comodule map with $0 \neq P$ of finite dimension. We know that $M \cong \oplus_{i \in I} E_H(S_i)$ for a set $\{S_i\}_{i \in I}$ of simple subcomodules of $M$
\cite[1.5h]{G}. There is $j \in I$ such that $g \vert_{E_H(S_j)}:E_H(S_j) \rightarrow P$ is nonzero. Composing the canonical projection $\pi_j:H \rightarrow E_H(S_j)$ with this map, we obtain that its image $N$ is a nonzero finite dimensional quotient comodule of $H$. By dualizing, $N^*$ is a finite dimensional left ideal of $H^*$. Then $0\neq N^* \subseteq \Rat(H^*)$ and by Theorem \ref{oldchar}, $H$ is co-Frobenius. \epf

\begin{remark}
The equivalence between (i) and (iv) is formulated in \cite[page 223]{Don1} for group schemes although its proof is completely different and strongly uses results and tools of group scheme theory. That $\Rat(H^*)\neq 0$ implies $H$ co-Frobenius is the key fact that allows us to prove this result in a much simpler manner. This must be seen as another instance of the power of the Fundamental Theorem of Hopf Modules.
\end{remark}

We are now in a position to prove our second main result:

\begin{theorem}\label{teo:cofrob}
The following assertions are equivalent:
\begin{enumerate}
\item[(i)] $H$ is co-Frobenius.
\item[(ii)] $H_{[0]}$ is co-Frobenius and $H$ is finitely generated as a left $H_{[0]}$-module.
\item[(iii)] $H_{[0]}$ is co-Frobenius and the standard filtration is finite.
\item[(iv)] The associated graded Hopf algebra $\gr H$ is co-Frobenius.
\item[(v)] $H_{[0]}$ is co-Frobenius and the diagram $R$ of $H$ is finite dimensional.
\end{enumerate}
Moreover, if $H_{[0]}$ is co-Frobenius, then $H_{[0]}=H_0 \stackrel{m}{...}H_0$ for some $m \geq 0$.
\end{theorem}

\pf (i) $\Rightarrow$ (ii) Since $H$ is co-Frobenius, its antipode is bijective \cite[Proposition 2]{R} and hence $\Ss(H_0)=H_0$. By Lemma \ref{lema:filt}, $H_{[0]}$ is a Hopf subalgebra of $H$; it is co-Frobenius because Hopf subalgebras inherit such a property \cite[Theorem 2.15]{Su2}. By Theorem \ref{teo:radford-beuatiful}, $H=H_0E_H(\ku) \subseteq H_{[0]}E_H(\ku) \subseteq H$.  \par \medskip

(ii) $\Rightarrow$ (iii) Let $h_1,...,h_r \in H$ be such that $H=H_{[0]}h_1+...+H_{[0]}h_r$. There is $m \geq 0$ such that $h_1,...,h_r \in H_{[m]}$. Then $H=H_{[0]}h_1+...+H_{[0]}h_r \subseteq H_{[0]}H_{[m]}=H_{[m]}$ by Lemma \ref{lema:filt}. \footnote{The hypothesis $\Ss(H_0) \subseteq H_0$ assumed in Lemma \ref{lema:filt} is not needed here. Nevertheless, $\Ss(H_0) \subseteq H_0$ holds. Since $H_{[0]}$ is co-Frobenius, $\Ss \vert_{H_{[0]}}$ is bijective and consequently $\Ss((H_{[0]})_{\stackrel{}{0}})=(H_{[0]})_{\stackrel{}{0}}$. Recall now that $(H_{[0]})_{\stackrel{}{0}}=H_0$.} \par \medskip

(iii) $\Rightarrow$ (i) Let $m \geq 0$ be minimal such that $H=H_{[m]}$. Since $H_{[0]}$ is co-Frobenius, the right $H_{[0]}$-comodule $H/H_{[m-1]}=H_{[m]}/H_{[m-1]}$ has a finite dimensional quotient $H_{[0]}$-comodule $M \neq 0$ by Theorem \ref{newchar1}. Then $M$ is a quotient $H$-comodule of $H$. Since $H$ is injective, Theorem \ref{newchar1} applies. \par \medskip

(iii) $\Rightarrow$ (iv) By hypothesis, $(\gr H)_{[0]}=H_{[0]}$ is co-Frobenius and the standard filtration of $H$ is finite. In view of Proposition \ref{stfilgrh} the standard filtration of $\gr H$ is finite. By (iii) $\Rightarrow$ (i), $\gr H$ is co-Frobenius. \par \medskip

(iv) $\Rightarrow$ (iii) Since $H_{[0]}$ is a Hopf subalgebra of $\gr H$, we have that $H_{[0]}$ is co-Frobenius. On the other hand, the standard filtration of $\gr H$ is finite by (i) $\Rightarrow$ (iii) applied to $\gr H$. From Proposition \ref{stfilgrh}, it follows that the standard filtration of $H$ is finite. \par \medskip

(iv) $\Leftrightarrow$ (v) The proof of this equivalence given in \cite[page 148]{AD} when $H_0$ is a subalgebra is also valid in this setting. The argument is as follows. For $r \in R$ write $\Delta_R(r)=r^{(1)} \otimes r^{(2)} \in R \otimes R$. Denoting the $H_{[0]}$-comodule structure map of $R$ by $\lambda:R \rightarrow H_{[0]} \otimes R$, set $\lambda(r)= r_{(-1)} \otimes r_{(0)}$. The comultiplication of $R \# H_{[0]}$ is given by $\Delta(r \# h)=(r^{(1)}\# {r^{(2)}}_{(-1)}h_{(1)}) \otimes ({r^{(2)}}_{(0)} \#  h_{(2)})$ for $r \# h\in R \# H_{[0]}$. Notice that if $K$ is a left coideal of $H_{[0]},$ then $R \# K$ is a left coideal of $R \# H_{[0]}$. \par \smallskip

By Proposition \ref{stfilgrh}, $(\gr H)_0=(H_{[0]})_{\stackrel{}{0}}=H_0$. Under the Hopf algebra isomorphism $\gr H \cong R \# H_{[0]},$ the coradical $(\gr H)_0$ corresponds to $R^0\# H_0=\ku \# H_0$.
In other words, there is a bijective correspondence between the set of isomorphism classes of simple $H$-comodules and the set of isomorphism classes of simple $\gr H$-comodules.
Take $\{S_i\}_{i \in I}$ a set of simple left coideals of $H$ such that $\gr H=\oplus_{i \in I} E_{\gr H}(S_i)$
and $H_{[0]}=\oplus_{i \in I} E_{H_{[0]}}(S_i)$. Then $R \# H_{[0]}=\oplus_{i \in I} R \# E_{H_{[0]}}(S_i)$ as left comodules.
This implies that $R \# E_{H_{[0]}}(S_i)$ is an injective left coideal of $R \# H_{[0]}$ containing $\ku \# S_i$.
Observe that $\ku \# S_i$ is the only simple left coideal contained in $R \# E_{H_{[0]}}(S_i)$: if $S \subset R \# E_{H_{[0]}}(S_i)$ is another one, then $S \subset (R \# E_{H_{[0]}}(S_i)) \cap (R \# H_{[0]})_0 \subset \ku \# S_i$, hence $S = \ku \# S_i$. Then  $E_{\gr H}(S_i)\cong R \# E_{H_{[0]}}(S_i)$. The claim now follows:

\noindent$\gr H$ is co-Frobenius $\iff \dim E_{\gr H}(\ku) < \infty \iff \dim R < \infty$ and $\dim E_{H_{[0]}}(\ku)< \infty \iff\dim R < \infty$ and $H_{[0]}$ is co-Frobenius.
\par \medskip

Finally, if $H_{[0]}$ is co-Frobenius, then $H_{[0]}=(H_{[0]})_{\stackrel{}{0}}E_{H_{[0]}}(\ku)$ by  Theorem \ref{teo:radford-beuatiful}. Recall from Lemma \ref{lema:filt} that $H_{[0]}=\bigcup_{r \geq 0} H_0^{(r)}.$ Since $\dim E_{H_{[0]}}(\ku)< \infty$, there is $t \geq 0$ such that $E_{H_{[0]}}(\ku) \subseteq H_0^{(t)}$. Then $H_{[0]}=(H_{[0]})_{\stackrel{}{0}}E_{H_{[0]}}(\ku)\subseteq H_0 H_0^{(t)}=H_0^{(t+1)}.$ \epf

\begin{remark}
Notice that our proof of (i) $\Leftrightarrow$ (iv) in Theorem \ref{teo:cofrob} is different to that in \cite{AD} when $H_0$ is a subalgebra. Given a simple left $H$-comodule $S$, a $\gr H$-comodule $\gr E_H(S)$ attached to $E_H(S)$ is constructed as $\gr E_H(S)=\oplus_{i \geq 0} \gr^i E_H(S)$ with $\gr^i E_H(S)=(E_H(S) \cap H_i)/(E_H(S) \cap H_{i-1})$. It is proved then there that $\gr E_H(S)\cong E_{{\gr H}}(S)$.
\end{remark}

\subsection{Exact sequences of co-Frobenius Hopf algebras}

In this subsection we will prove that the central term $B$ in an extension of Hopf algebras
$\ku \rightarrow A \rightarrow B \rightarrow C \rightarrow \ku$
is co-Frobenius if and only if $A$ and $C$ are co-Frobenius. We will also show that $B$ possesses a nonzero left integral $\int$  such that $\int \vert_{\stackrel{}{A}} \neq 0$ if and only if $A$ is co-Frobenius and $C$ is cosemisimple. We will derive from this that $B$ is cosemisimple if and only if $A$ and $C$ so are. \par \smallskip

The first result mentioned will be obtained as a consequence of Theorem \ref{newchar1} and another new characterization of co-Frobenius Hopf algebras that we present next. We will need the following description of the cotensor product.

\begin{lemma}\label{cotcoi} \cite[Lemma 3.1]{Sch1}
Let $M$ be a right $H$-comodule and $X$ be a left $H$-comodule. Let $X^{\bullet}$ denote $X$ but viewed as a right comodule via the antipode. Then
$M \Box_H X = (M \otimes X^{\bullet})^{\co H}$. \qed
\end{lemma}

\par \medskip

A homological condition characterizing co-Frobenius Hopf algebras is that the category of right (resp. left) comodules has enough projective objects, \cite[Theorems 3 and 10]{L}. In the particular case of the coordinate algebra of a group scheme, Donkin showed that the existence of a nonzero projective object suffices to characterize such Hopf algebras, \cite[Lemma 1]{Don2}. We observe that Donkin's result easily extends to arbitrary Hopf algebras.

\begin{theorem}\label{newchar2}
The following statements are equivalent:
\begin{enumerate}
\item[(i)] $H$ is co-Frobenius.
\item[(ii)] $\mh$ posseses a nonzero projective object.
\item[(iii)] Every injective right $H$-comodule is projective.
\end{enumerate}
\end{theorem}

\pf (i) $\Rightarrow$ (ii) By \cite{L}, as said above. \par \medskip

(ii) $\Rightarrow$ (iii) First we prove that every projective right $H$-comodule $M$ is injective. For $N,X,Y \in \mh$ with $X$ of finite dimension there is a natural isomorphism $\Hom_{H}(N \otimes X,Y) \cong \Hom_{H}(N, Y \otimes X^*),$ where $X^*$ is the left dual of $X$ constructed using the antipode. Then $N \otimes X$ is projective if $N$ is so. To show that $M$ is injective we must check that for an epimorphism $g: Z \rightarrow X$ of finite dimensional left $H$-comodules, the map $\id_M \Box_H g= (\id_M \otimes g) \vert_{M \Box_H Z}: M \Box_H Z \rightarrow M \Box_H X$ is surjective. It is known that the notions of \textit{injective} and \textit{coflat} are equivalent in the category of comodules over a coalgebra, \cite[Appendix, 2.1]{T3}. Since $M \otimes X^{\bullet}$ is projective, the map $\id_M \otimes g:M \otimes Z^{\bullet} \rightarrow M \otimes X^{\bullet}$ splits, so there exists a right $H$-comodule map $\theta:M \otimes X^{\bullet} \rightarrow M \otimes Z^{\bullet}$ such that $(\id_M \otimes g)\theta=\id_{M \otimes X^{\bullet}}$. Taking into account the inclusions $(\id_M \otimes g)((M \otimes Z^{\bullet})^{\co H}) \subseteq (M \otimes X^{\bullet})^{\co H}$ and $\theta((M \otimes X^{\bullet})^{\co H}) \subseteq (M \otimes Z^{\bullet})^{\co H}$  and Lemma \ref{cotcoi}, $\id \Box_H g$ splits and then it is surjective. Hence $M$ is injective. \par \smallskip

Let $P \in \mh$ be a nonzero projective and $Q$ a nonzero finite dimensional subcomodule of $P$. Consider the canonical map $\ku \rightarrow Q \otimes Q^*$. Let $S \in \mh$ be simple. We have an injective comodule map $S \rightarrow Q \otimes Q^* \otimes S \rightarrow P \otimes Q^* \otimes S$. The latter is projective, so it is injective by the previous paragraph. Then $E_H(S)$ is a direct summand of $P \otimes Q^* \otimes S$. Since $P \otimes Q^* \otimes S$ is projective, $E_H(S)$ is projective. Finally, every injective in $\mh$ is isomorphic to a direct sum of injective hulls of simple comodules, thus it is projective. \par \medskip

(iii) $\Rightarrow$ (i) By hypothesis, $E_H(\ku)$ is projective. It is known that a projective indecomposable comodule is finite dimensional, \cite[Lemma 1.2]{GN}.  \epf \vspace{-1pt}

We give an application of the previous theorem addressed to prove the announced result on exact sequences of Hopf algebras. Recall that a right $H$-comodule $M$ is said to be {\it finitely cogenerated} if there is a monomorphism of right $H$-comodules from $M$ into $H^n$ for some $n \in \Na$.

\begin{corollary}\label{cofquot}
Let $g:H \rightarrow K$ be a Hopf algebra map.
\begin{enumerate}
\item[(i)] If $H$ is co-Frobenius and $H$ is injective as right $K$-comodule, then $K$ is co-Frobenius.
\item[(ii)] If $H$ is finitely cogenerated as a right $K$-comodule and $K$ is co-Frobenius, then $H$ is co-Frobenius.
\end{enumerate}
\end{corollary}

\pf Let $\Res:\hm \rightarrow \km$ and $\Ind:=H \Box_K -:\km \rightarrow \hm$ be the restriction and induction functors respectively. We know that $\Res$ is left adjoint to $\Ind$. If $H$ is injective as a right $K$-comodule, then $\Ind$ is exact, and hence $\Res$ preserves projectives. On the other hand, $\Ind$ preserves injectives because $\Res$ is always exact. \par \smallskip

(i) By hypothesis and Theorem \ref{newchar2}, there is a nonzero projective object $P \in \hm$. Then $\Res P$ is a nonzero projective object in $\mk$ and, using again Theorem \ref{newchar2}, $K$ is co-Frobenius. \par \medskip

(ii) Let $f:H \rightarrow K^n$ be the monomorphism of $K$-comodules given by hypothesis. Take $M \in \km$ finite dimensional. We have a monomorphism of vector spaces $f \Box_K \id_M: H \Box_K M \rightarrow K^n\Box_K M \cong M^n$. From here, $\Ind M =H \Box_K M$ is finite dimensional. Since $K$ is co-Frobenius, $\Ind E_K(\ku)$ is a finite dimensional injective in $\hm$. Moreover, $\Ind(E_K(\ku)) \neq 0$ because it contains
$\Ind \ku = H \Box_K \ku =(H \otimes \ku^{\bullet})^{\co K} \cong H^{\co K} \neq 0$. By Theorem \ref{newchar2}, $H$ is co-Frobenius. \epf \vspace{4pt}

Recall from \cite{ADe} that a sequence of morphisms of Hopf algebras
$$\xymatrix{\ku \ar[r] & A \ar[r]^{\iota} & B \ar[r]^{\pi} & C \ar[r] & \ku}$$
is exact if $\iota$ is injective, $\pi$ is surjective,
\begin{align}
\label{eqn:exact1} \ker \pi &= BA^{+} \qquad \text{ and } \\
\label{eqn:exact2} B^{\co C}&= A.
\end{align}

There are some simplifications of this definition, see \cite{ADe, Sch2, T2}:

\begin{itemize}
  \item If $A$ is stable under the adjoint action of $B$ (i.e., $A$ is normal) and $B$ is faithfully flat as an $A$-module, then \eqref{eqn:exact1} implies \eqref{eqn:exact2}.

 \medbreak \item If $C$ is a quotient comodule of $B$ under the adjoint coaction (i.e., $C$ is conormal) and $B$ is faithfully coflat as a $C$-comodule, then \eqref{eqn:exact2} implies \eqref{eqn:exact1}.

 \medbreak   \item \emph{$A$ is a normal Hopf subalgebra of $B$ and $B$ is faithfully flat as an $A$-module} is equivalent to \emph{$B$ is faithfully coflat as a $C$-comodule and $C$ is a conormal quotient Hopf algebra of $B$.}

\end{itemize}

 \medbreak We are now ready to prove the announced result:

\begin{theorem}\label{exseq}
Let $\ku \rightarrow A \rightarrow B \rightarrow C \rightarrow \ku$ be an exact sequence of Hopf algebras with $B$ faithfully coflat as a $C$-comodule. Then, $B$ is co-Frobenius if and only if $A$ and $C$ are co-Frobenius.
\end{theorem}

\pf Assume that $A$ and $C$ are co-Frobenius. The functor $\Ind:=B \Box_C -:\mc \rightarrow \mb$ is exact because $B$ is injective as a right $C$-comodule. A nonzero left integral $\int^C:C \rightarrow \ku$ for $C$ is a surjective map of left $C$-comodules. Then $\Ind\int^C:\Ind C\rightarrow \Ind \ku$ in $\mb$ is surjective. Clearly $\Ind C \cong B$, and $\Ind \ku =B \Box_C \ku=(B \otimes \ku^{\bullet})^{\co C} \cong B^{\co C}=A$. So $A$ is a quotient of $B$ as a left $B$-comodule. Since $A$ is co-Frobenius, $A$ has a nonzero finite dimensional quotient left $A$-comodule (and hence $B$-comodule). Therefore $B$ has a nonzero finite dimensional quotient left $B$-comodule and by Theorem \ref{newchar1}, $B$ is co-Frobenius. \par \medskip

Conversely, if $B$ is co-Frobenius, then $A$ is co-Frobenius by \cite[Theorem 2.15]{Su2}. That $C$ is co-Frobenius follows from Corollary \ref{cofquot}(i) since $B$ is injective as a right $C$-comodule. \epf

\begin{examples}\label{examples}
(1) For commutative Hopf algebras, Sullivan proved Theorem \ref{exseq} by totally different methods in \cite[Theorem 2.20]{Su2}. \medbreak

(2) If the exact sequence $\ku \rightarrow A \rightarrow B \rightarrow C \rightarrow \ku$ is cleft, then $B$ is a bicrossproduct of $A$ and $C$. It is shown in \cite[Proposition 5.2]{BDGN} that $B$ is co-Frobenius if $A,C$ are so by checking that $\int^A \otimes \int^C$ is a nonzero integral for $B$ where $\int^A, \int^C$ are nonzero integrals for $A$ and $C$ respectively. \medbreak

In the next examples $\ku$ is an algebraically closed field of characteristic zero. \medbreak

(3) Let $G$ be a connected, simply connected, simple complex algebraic group and let $\epsilon$ be a primitive $\ell$-th root of 1, $\ell$
odd and $3 \nmid \ell$ if $G$ is of type $G_{2}$. It was shown in \cite[Example 4.1]{AD} using the Hopf socle that the quantum group $\Oe$ is co-Frobenius. An alternative proof follows from Theorem \ref{exseq}. For, $\Oe$ is noetherian and fits into an exact sequence
$\ku \rightarrow \Oc(G) \rightarrow \Oe \rightarrow  \bar{H} \rightarrow \ku$, where $\dim \bar{H} < \infty$ and $\Oc(G)$ is central in $\Oe$. Hence
$\Oe$ is faithfully flat over $\Oc(G)$ by \cite[Theorem 3.3]{Sch2}. \medbreak

(4) Let $\D= (I_{+}, I_{-}, N, \Gamma, \sigma, \delta)$ be a subgroup datum as in \cite[Definition 1.1]{AG} and let $A_{\D}$ be the Hopf algebra (quotient of $\Oe$) constructed in \cite[\S 2]{AG}. Then $A_{\D}$ is co-Frobenius iff the algebraic group $\Gamma$ is reductive. By \cite[Theorem 2.17]{AG}, $A_{\D}$ fits into the following diagram with exact rows and surjective vertical maps $$\xymatrix{\ku \ar[r] & \Oc(G) \ar[r]\ar@{>>}[d] & \Oe \ar@{>>}[d] \ar[r] & \bar{H}\ar@{>>}[d] \ar[r] & \ku\\
\ku \ar[r] &\Oc(\Gamma) \ar[r] &A_{\D} \ar[r] & H \ar[r] & \ku}$$
Hence $A_{\D}$ is noetherian, $\Oc(\Gamma)$ is central and $\dim H < \infty$; thus
$A_{\D}$ is faithfully flat over $\Oc(\Gamma)$ again by \cite[Theorem 3.3]{Sch2}. It is known that $\Oc(\Gamma)$ is co-Frobenius iff  $\Gamma$ is reductive \cite[Theorem 3]{Su2}. Then Theorem \ref{exseq} applies. \medbreak

(5) If $H$ is a co-Frobenius Hopf algebra and $\sigma: H \otimes H \to \ku$ is a convolution invertible 2-cocycle, then the twisted algebra $H^{\sigma}$ is again co-Frobenius, since the coalgebra structure remains unchanged, see \cite{DT} for details. In this way, many algebras of functions on multiparametric quantum groups are co-Frobenius, like those studied in \cite{AST}, which are twistings of $\Oc_{\epsilon}(GL(n))$. However, there are multiparametric quantum groups that do not arise as twistings as we point out next. Also, the twisting operation does not preserve quotient Hopf algebras. \medbreak

(6) Let $\ell$ be an odd natural number such that $\alpha^{-1}\beta$ is a primitive $\ell$-th root of unity and $\alpha^{\ell} = 1 = \beta^{\ell}$. The 2-parameter quantum group $\Oc_{\alpha,\beta}(GL(n))$ introduced in \cite{T} is co-Frobenius.
For, $\Oc_{\alpha,\beta}(GL(n))$ is noetherian and fits into an exact sequence
$\ku \rightarrow \Oc(GL(n)) \rightarrow \Oc_{\alpha,\beta}(GL(n)) \rightarrow  \bar{H} \rightarrow \ku$, where $\dim \bar{H} = \ell^{n^2}$ and $\Oc(GL(n))$ is central in $\Oc_{\alpha,\beta}(GL(n))$, see \cite[5.1 and 5.3]{Ga}. Then Theorem \ref{exseq} applies. Notice that these 2-parameter quantum groups are not twistings of the quantum $GL(n)$ discussed above, see [Ga, Remark 3.2(b)] and [T4, Theorem 2.6]. \medbreak

(7) Let $\D= (I_{+}, I_{-}, N, \Gamma, \sigma, \delta)$ be a subgroup datum as in \cite[Definition 1.1]{Ga} and let $A_{\D}$ be the Hopf algebra
constructed in \cite[Section 5.3]{Ga} (different to the mentioned above from \cite{AG}). Then $A_{\D}$ is co-Frobenius iff the algebraic group $\Gamma$ is reductive. By \cite[Theorem 5.23]{Ga}, $ A_{\D}$ fits into the following diagram with exact rows and surjective vertical maps
$$\xymatrix{\ku \ar[r] & \Oc(GL(n)) \ar[r]\ar@{>>}[d] & \Oc_{\alpha,\beta}(GL(n)) \ar@{>>}[d] \ar[r] & \bar{H}\ar@{>>}[d] \ar[r] & \ku\\
\ku \ar[r] &\Oc(\Gamma) \ar[r] &A_{\D} \ar[r] & H \ar[r] & \ku}$$
Then Theorem \ref{exseq} applies.
\end{examples}
\par \medskip

There is another result of Sullivan, \cite[Theorem 1.5]{Su2} (whose proof was apparently never published), stating that if $A \subset B$ are commutative Hopf algebras, then $B$ has a left (right) integral $\int$ such that $\int \vert_{\stackrel{}{A}} \neq 0$ if and only if $A$ has a nonzero integral and the quotient Hopf algebra $B/BA^+$ is cosemisimple. We next prove this theorem for exact sequences of arbitrary Hopf algebras. The following technical result will be needed:

\begin{lemma}
Let $H$ be a co-Frobenius Hopf algebra with nonzero left integral $\int$. Then:
\begin{enumerate}
\item[(i)] $E_H(\ku)\in \mh$ has a unique maximal subcomodule $M$.
\item[(ii)] $\int \vert_{\stackrel{}{E_H(\ku)}} \neq 0$ and $\int \vert_{\stackrel{}{M}}=0$.
\end{enumerate}
\end{lemma}

\pf (i) Let $M:=\Rad E_H(\ku)$ be the radical of $E_H(\ku)$, i.e., the intersection of all its maximal subcomodules.
Let $g \in H$ be the distinguished group-like element. By \cite[Theorem 5.2]{C}, $E_H(\ku)/\Rad E_H(\ku)\cong \ku g$. Since $\ku g$ is simple, $M$ is the unique maximal subcomodule of $E_H(\ku)$. \par \medskip

(ii) The distinguished group-like element satisfies $\int(h_{(1)})h_{(2)}=\int(h)g$ for all $h\in H$. Hence $f:H \rightarrow \ku g$, $h \mapsto \int(h)g$ is a morphism in $\mh$. Decompose $H=E_H(\ku)\oplus P$ as a right $H$-comodule. Then $\int(h_{(2)})h_{(1)}=\int(h)1_H \in E_H(\ku)\cap P$ for all $h \in P$. Hence $\int \vert_{\stackrel{}{P}}=0$. Since $\int \neq 0$, it must be $\int \vert_{\stackrel{}{E_H(\ku)}} \neq 0.$ Thus, $f\vert_{\stackrel{}{E_H(\ku)}}:E_H(\ku) \rightarrow \ku g$ is a nonzero morphism in $\mh$. Its kernel coincides with $M$ by (i) and so $\int \vert_{\stackrel{}{M}}=0$. \epf \vspace{3pt}

An extension of semisimple Hopf algebras is semisimple by \cite[Theorem 2.6(2)]{BM} and \cite[Theorem 2.2]{Sch3}, see also \cite[Proposition 3.1.18]{A-canadian}. Item (ii) of the next result extends this fact to cosemisimple Hopf algebras.

\begin{theorem}\label{exseq2}
Let $\ku \rightarrow A \rightarrow B \rightarrow C \rightarrow \ku$ be an exact sequence of Hopf algebras with $B$ faithfully coflat as a $C$-comodule.
\begin{enumerate}
\item[(i)] There is a nonzero left integral $\int$ for $B$ such that $\int \vert_{\stackrel{}{A}} \neq 0$ if and only if
$A$ is co-Frobenius and $C$ is cosemisimple.

\item[(ii)] The Hopf algebra $B$ is cosemisimple if and only if $A$ and $C$ so are.
\end{enumerate}
\end{theorem}

\pf Let notation be as in the proof of Theorem \ref{exseq}. \par \smallskip

(i) Assume that $A$ is co-Frobenius and $C$ is cosemisimple. We may choose $\int^C:C \rightarrow \ku$ splitting the inclusion map $i:\ku \rightarrow C$. Then the map of left $B$-comodules $\Ind \int^C:B  \rightarrow A$ splits $\Ind i:A \rightarrow B$. Under the previous isomorphisms $\Ind \ku \cong A$ and $\Ind C \cong B$, the map $\Ind i$ corresponds to the inclusion map of $A$ into $B$. So $A$ is isomorphic to a direct summand of $B$ as a left $B$-comodule. Set $B \cong A \oplus Q$ for some $Q \in \mb$. Since $A$ is co-Frobenius, there is a nonzero left integral $\int^A:A \rightarrow \ku$ (i.e., a map of left $A$-comodules). Then it is also a map of left $B$-comodules. The map $\int^B:B \rightarrow \ku$ defined by $\int^B\vert_{\stackrel{}{A}}=\int^A$ and $\int^B\vert_{\stackrel{}{Q}}=0$ is a nonzero left integral for $B$. \par \smallskip

Conversely, it is clear that $A$ is co-Frobenius. We prove that $C$ is cosemi\-simple.
We can take the injective hull $E_A(\ku)$ as a $B$-subcomodule of $E_B(\ku)$ (viewing $E_A(\ku)$ as a $B$-comodule). Suppose that $E_A(\ku) \neq E_B(\ku)$. Then $E_A(\ku) \subseteq M$, with $M$ the unique maximal subcomodule of $E_B(\ku)$. From the hypothesis and the precedent lemma, we have $0 \neq \int \vert_{\stackrel{}{A}}(E_A(\ku))=\int(E_A(\ku))\subseteq \int(M)=0,$ a contradiction. Therefore $E_A(\ku)=E_B(\ku)$. This means that $E_A(\ku)$ is injective when viewed as a right $B$-comodule. If $S \in \ma$ is simple, we know that $E_A(S)$ is a direct summand of $S \otimes E_A(\ku)$ as an $A$-comodule (hence as a $B$-comodule either). Since the latter is injective as a $B$-comodule, $E_A(S)$  so is. This implies that $A$ is injective as a right $B$-comodule. There is a right $B$-comodule $Q$ such that $B \cong A \oplus Q$.
Applying the restriction functor $\Res:\mbr \rightarrow \mcr$, we get  $\Res B \cong \Res A \oplus \Res Q$. Taking into account that $B^{\co C}=A$, we have that $\Res A$ is isomorphic to a direct sum of copies of $\ku$. As $\Res B$ is injective, from the above, $\ku$ is injective as a $C$-comodule, and so $C$ is cosemisimple. \par \medskip

(ii) If $B$ is cosemisimple, a nonzero left integral $\int$ for $B$ satisfies $\int(1_B) \neq 0$. Then $\int \vert_{\stackrel{}{A}}(1_A) \neq 0$, giving that $A$ and $C$ are cosemisimple. Finally, if $A$ and $C$ are cosemisimple, by (i), there exists a left integral $\int$ for $B$ such that $\int \vert_{\stackrel{}{A}} \neq 0$. Since $A$ is cosemisimple, $0 \neq \int \vert_{\stackrel{}{A}}(1_A)=\int(1_B)$. From this, $B$ is cosemisimple. \epf

\begin{remark}
Notice that the hypothesis of $B$ being faithfully coflat as a $C$-comodule was not used in the proof of the implication from right to left in both statements.
\end{remark}

\subsection{Finite dual co-Frobenius Hopf algebras}
In this last subsection we give one more application of Theorem \ref{newchar2}. We obtain a result dual to Corollary \ref{cofquot} for finite dual Hopf algebras. We previously need the dual version of Lemma \ref{cotcoi} that appears in \cite[Lemma 4.1]{Sch1}. \par \medskip

Let $K$ be a Hopf algebra. Given a right $K$-module $M$ denote the quotient vector space $M/MK^+$ by $\overline{M}$.

\begin{lemma}\label{teninv}
Let $M$ and $X$ be right and left $K$-modules respectively. Let $X^{\bullet}$ denote $X$ but viewed as a right module via the antipode. Then $\overline{M \otimes X^{\bullet}} \cong M \otimes_{K} X$. \qed
\end{lemma}

As usual $H^0$ denotes the finite or Sweedler dual of $H$, i.e., the subspace of $H^*$ spanned by the matrix coefficients of all finite dimensional $H$-modules.

\begin{proposition}\label{cofindual}
Let $g:K \rightarrow H$ be a Hopf algebra map. Assume that $H$ is finitely generated as a right $K$-module. Then
\begin{enumerate}
\item[(i)] If $K^0$ is co-Frobenius, then $H^0$ is co-Frobenius.
\item[(ii)] If $H^0$ is co-Frobenius and $H$ is flat as a right $K$-module, then $K^0$ is co-Frobenius.
\end{enumerate}

\end{proposition}

\pf Let $\mkfd$ and $\mhfd$ denote the categories of finite dimensional left $K$-modules and $H$-modules respectively. We may identify
$\mkfd$ as the full subcategory of finite dimensional objects in $\kcerom$.

Since $H_K$ is finitely generated, there is an epimorphism of right $K$-modules $g:K^n \rightarrow H$ for some $n \in \Na$.
If $M \in \mkfd$, then the induced linear map $g  \otimes_K \id_M: M^n\cong  K^n \otimes_K M \rightarrow H \otimes_K M$ is surjective and hence $H \otimes_K M$ is finite dimensional. Thus we may consider the induction and restriction functors $\Ind=H \otimes_K - \hspace{-2pt}:\hspace{-3pt}\mkfd \rightarrow \mhfd$ and $\Res\hspace{-2pt}:\hspace{-3pt}\mhfd \rightarrow \mkfd$. We know that $\Res$ is right adjoint to $\Ind$. The functor $\Ind$ preserves projectives
because $\Res$ is exact. \par \smallskip

(i) The projective cover $P(\ku)$ of $\ku$ in $\kcerom$ belongs to $\mkfd$, see \cite[Theorem 10, Lemma 15]{L}. Then $\Ind P(\ku) \in \mhfd$ is projective. It will be nonzero if we show that its quotient $\Ind\ku$ is nonzero. For, we apply Lemma \ref{teninv} to obtain $\Ind \ku =H \otimes_K \ku \cong \overline{H \otimes \ku^{\bullet}} \cong \overline{H}=H/HK^+$. The latter is nonzero since $1_H \notin HK^+$. The projectives in $\mhfd$ coincide with the finite dimensional projectives in $\hcerom$ (this follows from the local finiteness of comodules). By Theorem \ref{newchar2}, $H^0$ is co-Frobenius. \par \smallskip

(ii) As $H$ is flat as a right $K$-module, $\Ind$ is exact and so $\Res$ preserves injectives. If $H^0$ is co-Frobenius, there is a nonzero finite dimensional injective $Q \in \hcerom$ by Theorem \ref{newchar2}. So $\Res Q$ is injective in $\mkfd$. Taking into account that the injectives in $\mkfd$ are exactly the finite dimensional injectives in $\kcerom$, and using once again Theorem \ref{newchar2}, $K^0$ is co-Frobenius. \epf

Examples of co-Frobenius (indeed cosemisimple) Hopf algebras were constructed in \cite[Corollary 3.3]{C} as finite dual Hopf algebras of group algebras of locally finite groups whose elements have order not divisible by $char(\ku)$. More generally, it was shown there that if $H$ is a Hopf algebra that is Von Neumann regular as an algebra, then $H^0$ is cosemisimple. \par \medskip

\bigskip
\begin{center}
{\bf ACKNOWLEDGMENT}
\end{center}

N. A. was partially supported by ANPCyT-Foncyt, CONICET, Ministerio de Ciencia y
Tecnolog\'{\i}a (C\'ordoba) and Secyt-UNC. J. C. is supported by projects MTM2008-03339 from MICINN and FEDER and P07-FQM03128 from Junta de Andaluc\'{\i}a. This work was carried out in the framework of the Spanish American joint project A/4742/06 from AECID. We thank Vyjayanthi Chari and Jacob Greenstein for interesting discussions on Hopf algebras related to quantum affine algebras.

\end{document}